\documentclass[reqno,12pt]{amsart}
\usepackage{amscd}
\usepackage[mathscr]{eucal}
\usepackage{amssymb}
\usepackage{latexsym}
\usepackage{amsthm}
\usepackage{amsmath}
\usepackage{graphicx}
\usepackage{pb-diagram}

\theoremstyle{plain}

\newtheorem{thm}{Theorem}[section]

\newtheorem{cor}{Corollary}[section]
\newtheorem{lem}{Lemma}[section]

\theoremstyle{definition}

\newtheorem{dfn}{Definition}[section]
\newtheorem{prf}{Proof}[section]

\DeclareMathOperator*{\esssup}{ess\, sup}

\title{On almost convergence on the real line and its application to bounded analytic functions}
\author{Ryoichi Kunisada}

\address{}

\email{tk-waseda@ruri.waseda.jp}
\subjclass[2010]{}
\date{}

\begin{document}

\begin{abstract}
We address the study of topologically invariant means and almost convergence on the real numbers $\mathbb{R}$. Here, the former is a certain class of invariant means on $L^{\infty}(\mathbb{R})$ and the latter is a summability method defined by them. Almost convergence on $\mathbb{R}$ was firstly introduced by Raimi (1957) as a generalization of Lorentz's almost convergence for bounded sequences. We extensively generalize his result of analytic characterization of almost convergence and explore its application to the theory of Hardy space. Specifically, we establish the relation between the asymptotic behavior on the imaginary axis and that at infinity of  bounded analytic functions defined on the right half plane.
\end{abstract}
\maketitle

\bigskip
\section{Introduction}
We study a certain summability method which we call almost convergence. This notion was firstly introduced by Lorentz for bounded functions on the additive semigroup of nonnegative integers $\mathbb{Z}_+$ (see [9]). After that, several authors generalized this notion to general locally compact amenable groups or semigroups (see [1], [2], [13], [14]). Here, we study exclusively almost convergence for the additive group of real numbers $\mathbb{R}$. This can be viewed as a continuous version of Lorentz's almost convergence and essentially equivalent to the one introduced by Raimi in [13]. In the paper, he gave a necessary and sufficient condition for essentially bounded measurable functions on $\mathbb{R}$ to be almost convergent, which is analogous to the one given by Lorentz. One of the main objectives of this paper is to generalize Raimi's result to obtain a more general form of necessary and sufficient condition for almost convergence including his result as a special case. Furthermore, we provide an application of almost convergence to the study of the Hardy space of bounded analytic functions on the right half plane. 

The paper organized as follows. In Section 2, we give some definitions and preliminary results which are needed in the later sections. In Section 3, we study topologically invariant means in detail and provide a necessary and sufficient condition for almost convergence, which is one of the main result of  this paper. This is very general and we can obtain many of analytic conditions including the known result due to Raimi. Our argument is based on the fact that topologically invariance is characterized with invariance and a kind of Fubini type property.

Section 4 deals with the Hardy space $H^{\infty}(\mathbb{C}^+)$ of bounded analytic functions on the right half plane. Using the results in Section 3, we see that the behaviour of a function in $H^{\infty}(\mathbb{C}^+)$ at infinity can be expressed in terms of almost convergence of its boundary function. We also treat the relation between topologically invariant means on $L^{\infty}(\mathbb{R})$ and the maximal ideal space of $H^{\infty}(\mathbb{C}^+)$. 

\section{Preliminaries}
Let $L^1(\mathbb{R})$ be the group algebra of $\mathbb{R}$, $L^{\infty}(\mathbb{R})$ be the set of essentially bounded functions on $\mathbb{R}$ and $C_u(\mathbb{R})$ be the set of bounded, uniformly continuous functions on $\mathbb{R}$.
Let us denote a general element of $L^1(\mathbb{R})$ by the symbols $f, g, \cdots$ and that of $L^{\infty}(\mathbb{R})$(and $C_u(\mathbb{R})$) by $\phi, \psi, \cdots$. For functions $f \in L^1(\mathbb{R})$ and $\phi \in L^{\infty}(\mathbb{R})$, we denote as $f_s(x) := f(x+s)$ and $\phi_s(x) := \phi(x+s)$, the translates of $f$ and $\phi$ by $s \in \mathbb{R}$, respectively. Let $\mathcal{M}(\mathbb{R}) \; (\mathcal{M}_0(\mathbb{R}))$ be the set of means on $L^{\infty}(\mathbb{R}) \; (C_u(\mathbb{R}))$, that is, the elements $\Phi$ of the dual space $L^{\infty}(\mathbb{R})^* \; (C_u(\mathbb{R})^*)$ such that  

\smallskip
\noindent
$\mathrm{(i)}$ $\Phi$ is positive, i.e., $\Phi(\phi) \ge 0$ whenever $\phi \ge 0$; \\
$\mathrm{(ii)}$ $\varphi(1) = 1$, where $1$ is the constant function that takes the value $1$ everywhere.  

\smallskip
\noindent
A mean $\Phi$ on $L^{\infty}(\mathbb{R}) \; (C_u(\mathbb{R}))$ is said to be invariant if it satisfies 

\noindent
\smallskip
$\mathrm{(iii)}$ $\Phi(\phi_s) = \Phi(\phi)$ for every $s \in \mathbb{R}$. 

\smallskip
\noindent
The set of invariant means on $L^{\infty}(\mathbb{R}) \; (C_u(\mathbb{R}))$ is denoted by $\mathcal{I}(\mathbb{R}) \; (\mathcal{I}_0(\mathbb{R}))$. Now we define topologically invariant means, which is the main objective of this paper. Let $P(\mathbb{R})$ be the subset of $L^1(\mathbb{R})$ consisting of those elements such that $f \ge 0$ and $\int_{\mathbb{R}} f(x)dx = 1$. A mean $\Phi$ on $L^{\infty}(\mathbb{R})$ is said to be topologically invariant if it satisfies the following condition ([5]): 

\smallskip
\noindent
$\mathrm{(iv)}$ $\Phi(f * \phi) = \Phi(\phi)$ for every $f \in P(\mathbb{R})$ and $\phi \in L^{\infty}(\mathbb{R})$, where $f * \phi$ is the convolution of $f$ and $\phi$ defined by
\[f * \phi(x) = \int_{\mathbb{R}} \phi(x-t)f(t)dt \quad (x \in \mathbb{R}). \]
Note that $f * \phi$ is in $C_u(\mathbb{R})$. Let us denote the set of all topologically invariant means on $L^{\infty}(\mathbb{R})$ by $\mathcal{T}(\mathbb{R})$. It is easy to see $\mathcal{T}(\mathbb{R}) \subseteq \mathcal{I}(\mathbb{R})$, in fact, if $\Phi \in \mathcal{T}(\mathbb{R})$, we have
\[\Phi(\phi_s) = \Phi(f * \phi_s) = \Phi(f_s * \phi) = \Phi(\phi) \]
and thus $\Phi \in \mathcal{I}(\mathbb{R})$. For more detailed account of invariant and topologically invariant means, see [3], [12].
We now define almost convergence for functions in $L^{\infty}(\mathbb{R})$ as follows.

\begin{dfn}
Let $\phi$ be in $L^{\infty}(\mathbb{R})$. $\phi$ is almost convergent to a complex number $\alpha$ if 
\[\Phi(\phi) = \alpha \]
holds for every $\Phi \in \mathcal{T}(\mathbb{R})$. In this case, we write as $\phi \xrightarrow{ac} \alpha$.
\end{dfn}

Let $L^{\infty}_R(\mathbb{R})$ be the set of real-valued essentially bounded functions on $\mathbb{R}$. Writing $\phi = u + iv$, where $u, v \in L^{\infty}_R(\mathbb{R})$, we obviously have $\phi \xrightarrow{ac} \alpha + i\beta$ if and only if $u \xrightarrow{ac} \alpha \ \text{and} \ v \xrightarrow{ac} \beta$.

In the following section, we provide a necessary and sufficient condition for a given mean on $L^{\infty}(\mathbb{R})$ to be topologically invariant. For our purpose, the following result which is derived from the Hahn-Banach theorem will be an essential tool.
\begin{thm}
Let $X$ be a real locally convex space. Let $\mathcal{C}$ be a compact convex subset of $X$ and $S$ is a subset of $\mathcal{C}$. Then,  the closed convex hull $\overline{co}(S)$ of $S$ is equal to $\mathcal{C}$ if and only if
\[\sup_{x \in S} \varphi(x) = \sup_{x \in \mathcal{C}} \varphi(x) \]
for every $\varphi \in X^*$, the dual space of $X$.
\end{thm}

\begin{prf}
Necessity is obvious by the linearlity and the continuity of $\varphi \in X^*$. Assume that $\overline{co}(S) \subsetneq \mathcal{C}$. Take $x_0 \in \mathcal{C} \setminus \overline{co}(S)$. Then, by the Hahn-Banach theorem, there exists a $\varphi \in X^*$ such that 
\[\sup_{x \in \overline{co}(S)} \varphi(x) < \varphi(x_0). \]
This implies that 
\[\sup_{x \in S} \varphi(x) = \sup_{x \in \overline{co}(S)} \varphi(x) < \varphi(x_0) \le \sup_{x \in \mathcal{C}} \varphi(x), \]
completing the proof.
\end{prf}

\begin{cor}
Let $\mathcal{C}$ be a weak* compact convex subset of $\mathcal{M}$ and let $S$ be a subset of $\mathcal{C}$. Then, $\mathcal{C}=\overline{co}(S)$ if and only if 
\[\sup_{\Phi \in S} \Phi(\phi) = \sup_{\Phi \in \mathcal{C}} \Phi(\phi) \]
holds for every $\phi \in L^{\infty}_R(\mathbb{R})$.
\end{cor}

\begin{prf}
Since necessity is obvious, we show sufficiency. For each element $\Phi$ of $\mathcal{C}$, let $\Phi_R$ be the restriction of $\Phi$ to the real Banach space $L^{\infty}_R(\mathbb{R})$. Set $\mathcal{C}_R = \{\Phi_R : \Phi \in \mathcal{C}\}$ and $S_R = \{\Phi_R : \Phi \in S\}$. Then, by the assumption, we can apply Theorem $2.1$ to $\mathcal{C}_R$ and $S_R$ and obtain $\overline{co}(S_R) = \mathcal{C}_R$.  It means that for each $\Phi \in \mathcal{C}$, there exists a net $\{\Phi_{\alpha}\}$ in $co(S)$ such that
\[\lim_{\alpha} (\Phi_{\alpha})(\phi) = \Phi(\phi) \]
holds for every $\phi \in L^{\infty}_R(\mathbb{R})$. From this equation, we have
\[\lim_{\alpha} \Phi_{\alpha}(\phi) = \lim_{\alpha} \Phi_{\alpha}(u+iv) = \lim_{\alpha} \{\Phi_{\alpha}(u) + i\Phi_{\alpha}(v)\} = \Phi(u) + i\Phi(v) = \Phi(\phi) \]
for every $\phi = u + iv \in L^{\infty}(\mathbb{R})$. This means that $w^*\mathchar`-\lim_{\alpha} \Phi_{\alpha} = \Phi$ and it implies that $\overline{co}(S) = \mathcal{C}$. This proves the theorem.
\end{prf}

We remark a fact which is a direct consequence of Corollary 2.1. Let $\mathfrak{C}$ be the set of compact convex subsets of $\mathcal{M}$ and $\mathfrak{P}$ be the set of sublinear functionals $\overline{p}$ on $L^{\infty}_{\mathbb{R}}(\mathbb{R})$ such that $\phi \ge 0$ implies $\overline{p}(\phi) \ge 0$ and $\overline{p}(1) = 1$. Then, note that the two partically ordered sets $(\mathfrak{C}, \subset)$ and $(\mathfrak{P}, \le)$ are isomorphic via the following correspondence:
\[\mathfrak{C} \ni \mathcal{C} \mapsto \overline{p}(\phi) = \sup_{\Phi \in \mathcal{C}} \Phi(\phi) \in \mathfrak{P}, \]
\[\mathfrak{P} \ni \overline{p} \mapsto \mathcal{C} = \{\Phi \in \mathcal{M} : \Phi(\phi) \le \overline{p}(\phi) \ (\forall \phi \in L^{\infty}_R(\mathbb{R}))\}. \]

The following elementary lemma plays an important role.
\begin{lem}
Let $\Phi \in \mathcal{M}$. For any $\phi \in L^{\infty}_R(\mathbb{R})$, we have
\[\Phi(\phi) \le \esssup_{x \in \mathbb{R}} \phi(x). \]
\end{lem}

\begin{prf}
Put $\alpha = \esssup_{x \in \mathbb{R}} \phi(x)$. Then, by the properties $\mathrm{(i), (ii)}$ of means and the fact that $\alpha - \phi \ge 0$, we obtain
\[\Phi(\alpha-\phi) \ge 0 \Leftrightarrow \alpha = \Phi(\alpha) \ge \Phi(\phi), \]
completing the proof.
\end{prf}

\section{Topologically invariant means}
The importance of the topologically invariant means comes from the following property: 
\begin{equation}
\Phi(f * \phi) = \int_{\mathbb{R}} \Phi(\phi_{-t})f(t)dt,    \tag{\dag} 
\end{equation}
where $\phi \in L^{\infty}(\mathbb{R})$, $f \in L^1(\mathbb{R})$ and $\Phi \in L^{\infty}(\mathbb{R})^*$, the dual space of $L^{\infty}(\mathbb{R})$. 

Note that the topologically invariant means can be characterized as invariant means with the property $(\dag)$. 
In fact, neccesity is obvious by the fact that $L^1(\mathbb{R})$ is spaned by $P(\mathbb{R})$. Conversely, if $\Phi \in \mathcal{M}$ satisfies the invariance and the property $(\dag)$, then, for any $f \in P(\mathbb{R})$, we have
\[\Phi(f * \phi) = \int_{\mathbb{R}} \Phi(\phi_{-t})f(t)dt = \int_{\mathbb{R}} \Phi(\phi)f(t)dt = \Phi(\phi)\int_{\mathbb{R}} f(t)dt = \Phi(\phi), \]
which shows that $\Phi$ is topologically invariant.

The property $(\dag)$ does not hold in general, but it is always true for $C_u(\mathbb{R})^*$. 
\begin{lem}
Let $f \in L^1(\mathbb{R}), \phi \in C_u(\mathbb{R})$ and $\Phi \in C_u(\mathbb{R})^*$. Then, we have
\[\Phi(f * \phi) = \int_{\mathbb{R}} \Phi(\phi_{-t})f(t)dt. \]
\end{lem}

\begin{prf}
Since the mapping $\mathbb{R} \ni s \mapsto \phi_s \in C_u(\mathbb{R})$ is continuous, the result follows from the well known fact in Bochner integral theory (see [19]).
\end{prf}
There exists a special class of means on $L^{\infty}(\mathbb{R})$ satisfying the property $(\dag)$. Let $h \in P(\mathbb{R})$ and define $\Phi_h \in \mathcal{M}$ by
\[\Phi_h(\phi) = \int_{\mathbb{R}} \phi(t)h(-t)dt, \]
where $\phi \in L^{\infty}(\mathbb{R})$. Then, the validity of the property $(\dag)$ follows from Fubini's theorem. Below, we introduce an extended class of means on $L^{\infty}(\mathbb{R})$ satisfying the property $(\dag)$, which contains the above examples.

For each $f \in P(\mathbb{R})$, define the functionals $\overline{F}$ and $\underline{F}$ on $L^{\infty}_R(\mathbb{R})$ by
\[
\overline{F}(\phi) = \sup_{x \in \mathbb{R}} (f * \phi)(x), 
\]
\[
\underline{F}(\phi) = \inf_{x \in \mathbb{R}} (f * \phi)(x), 
\]
where $\phi \in L^{\infty}_R(\mathbb{R})$. Note that $\overline{F}$ is sublinear and the relation $\underline{F}(\phi) = -\overline{F}(-\phi)$ holds.

Let $\mathcal{F}$ be a weak* compact convex subset of $\mathcal{M}$ consisting of those elements $\Phi$ such that
\[\Phi(\phi) \le \overline{F}(\phi) \]
for every $\phi \in L^{\infty}_R(\mathbb{R})$. 
\begin{thm}
For a mean $\Phi$ on $L^{\infty}(\mathbb{R})$, $\Phi \in \mathcal{F}$ if and only if there exists a mean $\Phi_0$ on $C_u(\mathbb{R})$ such that 
\[\Phi(\phi) = \Phi_0(f * \phi) \]
for every $\phi \in L^{\infty}(\mathbb{R})$.
\end{thm}

\begin{prf}
By Corollary $2.1$, it is sufficient to show that 
\[\sup_{\Phi_0 \in \mathcal{M}_0} \Phi_0(f * \phi) = \overline{F}(\phi) \]
for every $\phi \in L^{\infty}(\mathbb{R})$. 
By Lemma $2.1$, for any $\Phi_0 \in \mathcal{M}_0$, we have
\[\Phi_0(f * \phi) \le \esssup_{x \in \mathbb{R}} (f * \phi)(x) = \sup_{x \in \mathbb{R}} (f * \phi)(x) = \overline{F}(\phi). \]
Next, we show the reverse inequality. Put $\alpha = \esssup_{x \in \mathbb{R}} f * \phi(x) = \overline{F}(\phi)$. Then, there is a sequence $\{x_n\}$ in $\mathbb{R}$ such that 
\[\lim_{n \to \infty} (f * \phi)(x_n) = \alpha. \]
Consider the sequence $\{\hat{x}_n\}$ of evalution mappings at $x_n$, that is, $\hat{x}_n(\phi) := \phi(x_n)$. Let $\Phi^{\prime}$ be a cluster point of the sequence $\{\hat{x}_n\}$ in $C_u(\mathbb{R})^*$. Then, obviously, we have
\[\Phi^{\prime}(f*\phi) = \alpha, \]
which means that 
\[\sup_{\Phi_0 \in \mathcal{M}_0} \Phi_0(f * \phi) \ge \overline{F}(\phi). \]
We complete the proof.
\end{prf}

\begin{thm}
For any $f \in P(\mathbb{R})$, each element of $\mathcal{F}$ satisfies the property $(\dag)$.
\end{thm}

\begin{prf}
By Theorem $3.1$, for any $\Phi \in \mathcal{F}$, there exists a $\Phi_0 \in \mathcal{M}_0$ such that 
\[\Phi(\phi) = \Phi_0(f * \phi), \]
where $\phi \in L^{\infty}(\mathbb{R})$. Then, for any $g \in L^1(\mathbb{R})$, by Lemma $3.1$, we have
\begin{align}
\Phi(g * \phi) &= \Phi_0(f * g * \phi) = \Phi_0(g * f * \phi) \notag \\
&= \int_{\mathbb{R}} \Phi_0((f * \phi)_{-t}))g(t)dt = \int_{\mathbb{R}} \Phi_0(f * \phi_{-t})g(t)dt \notag \\
&= \int_{\mathbb{R}} \Phi(\phi_{-t})g(t)dt. \notag 
\end{align}
We obtain the result.
\end{prf}

\begin{thm}
Let $f$ be in $P(\mathbb{R})$. For any $\Phi \in \mathcal{M}$, $\Phi \in \mathcal{T}$ if and only if $\Phi$ is invariant and $\Phi \in \mathcal{F}$.
\end{thm}

\begin{prf}
Sufficiency is obvious from the fact mentioned above and Theorem $3.2$. We show necessity. Let $\Phi$ be a topologically invariant mean. Then, the invariance of $\Phi$ is obvious. Since $\Phi$ is a mean, we have
\begin{align}
\Phi(\phi) &\le \esssup_{x \in \mathbb{R}} \phi(x) \notag 
\end{align}
where $\phi \in L^{\infty}_R(\mathbb{R})$. Since $\Phi$ is topologically invariant, $\Phi(\phi) = \Phi(f * \phi)$ holds and thus, we have
\[\Phi(\phi) = \Phi(f * \phi) \le \sup_{x \in \mathbb{R}} (f * \phi)(x) = \overline{F}(\phi). \]
Hence $\Phi$ is in $\mathcal{F}$ and we complete the proof.
\end{prf}

Let $f$ be in $P(\mathbb{R})$. For each $r > 0$, we define $f_r(x) = r^{-1}f(x/r)$. Let $\overline{F}_u $ and $\underline{F}_u$ be the functionals on $L^{\infty}_{\mathbb{R}}(\mathbb{R})$ defined by
\[\overline{F}_u(\phi) = \limsup_{r \to \infty} \sup_{x \in \mathbb{R}} (f_r * \phi)(x) = \limsup_{r \to \infty} \overline{F}_r(\phi), \]
\[\overline{F}_u(\phi) = \liminf_{r \to \infty} \inf_{x \in \mathbb{R}} (f_r * \phi)(x) = \liminf_{r \to \infty} \underline{F}_r(\phi), \]
where we set $\overline{F_r}(\phi) = \sup_{x \in \mathbb{R}} (f_r * \phi)(x)$ and $\underline{F}_r(\phi) = \inf_{x \in \mathbb{R}} (f_r * \phi)(x)$ for each $r > 0$. Then, note that $\overline{F}_u$ is sublinear and the relation $\underline{F}_u(\phi) = -\overline{F}_u(-\phi)$ holds.

The most interesting examples of such functionals are the cases $D(x) = I_{[-1,1]}(x)$, the characteristic function of the interval $[-1, 1]$, and $P(x) = \frac{1}{\pi}\frac{1}{1+x^2}$. Note that the functions $P_x(y) = \frac{1}{\pi}\frac{1}{x}\frac{1}{1+\left(\frac{y}{x}\right)^2} = \frac{1}{\pi}\frac{x}{x^2 +y^2} (x > 0, y \in \mathbb{R})$ is the Poisson kernel.
These kernels give the following sublinear functionals on $L^{\infty}_R(\mathbb{R})$:
\begin{align}
\overline{D_u}(\phi) &= \limsup_{\theta \to \infty} \sup_{x \in \mathbb{R}} \int_{-\infty}^{\infty} \phi(x-t)D_{\theta}(t)dt \notag \\
&= \limsup_{\theta \to \infty} \sup_{x \in \mathbb{R}} \frac{1}{2\theta}\int_{x-\theta}^{x+\theta} \phi(t)dt, \notag
\end{align}

\begin{align}
\overline{P_u}(\phi) &= \limsup_{x \to \infty} \sup_{y \in \mathbb{R}} \int_{-\infty}^{\infty} \phi(y-t)P_x(t)dt \notag \\
&= \limsup_{x \to \infty} \sup_{y \in \mathbb{R}} \frac{1}{\pi}\int_{-\infty}^{\infty} \phi(t)\frac{x}{x^2+(y-t)^2}dt. \notag 
\end{align}
Here we mention the fact that $\limsup$ can be replaced with $\lim$ in the above two formulas. The following theorem is the main result of this section.
\begin{thm}
Let $f \in P(\mathbb{R})$ and $\Phi \in \mathcal{M}$. Then, $\Phi$ is topologically invariant if and only if 
\[\Phi(\phi) \le \overline{F_u}(\phi) \]
for every $\phi \in L^{\infty}_R(\mathbb{R})$.
\end{thm}

\begin{prf}
First, we show sufficiency. Suppose that for any $\phi \in L^{\infty}_R(\mathbb{R})$, 
\[\Phi(\phi) \le \overline{F}_u(\phi) \]
holds true. Since $\overline{F}_u \le \overline{F}$ is valid by definition, Then, by Theorem 3.2, $\Phi$ satisfies the property $(\dag)$. Thus, it remains to show that $\Phi$ is invariant. To this end, it is sufficient to show that 
\[\lim_{r \to \infty} \int_{-\infty}^{\infty} |f_r(t) - f_r(t+s)|dt = 0 \]
for every $s \in \mathbb{R}$. In fact, if this holds true, we have 
\begin{align}
\Phi(\phi-\phi_s) &\le \overline{F}_u(\phi-\phi_s) \notag \\
&= \limsup_{r \to \infty} \sup_{x \in \mathbb{R}} \{(f_r * \phi)(x) - (f_r * \phi)(x+s)\} \notag \\
&= \limsup_{r \to \infty} \sup_{x \in \mathbb{R}} \int_{-\infty}^{\infty} \phi(x-t)\{f_r(t)-f_r(t+s)\}dt \notag \\
&\le \limsup_{r \to \infty} \|\phi\|_{\infty} \int_{-\infty}^{\infty} |f_r(t)-f_r(t+s)|dt = 0. \notag 
\end{align}
In the same way, we have
\[\Phi(\phi-\phi_s) \ge \underline{F}_u(\phi-\phi_s) = -\overline{F}_u(\phi_s-\phi) = 0, \]
which inequalities show that $\Phi(\phi-\phi_s) = 0$ and we obtain the desired result. Now we prove the above equation. By substitution of variable, we obtain
\begin{align}
\int_{-\infty}^{\infty} |f_r(t) - f_r(t+s)|dt &= \int_{-\infty}^{\infty} \frac{1}{r}\left|f\left(\frac{t}{r}\right)-f\left(\frac{t+s}{r}\right)\right|dt \notag \\
&= \int_{-\infty}^{\infty} |f(t) - f(t+s/r)|dt. \notag 
\end{align}
The last equation tends to $0$ as $r$ tends to infinity for any fixed $s \in \mathbb{R}$. We obtain the result.

Next, we show necessity. Suppose $\Phi$ is topologically invariant. Then, by the same argument as in the proof of Theorem 3.3, we obtain
\[\Phi(\phi) \le \overline{F}_r(\phi) \]
for every $\phi \in L^{\infty}_R(\mathbb{R})$ and $r > 0$. Then, taking the limit superior as $r \to \infty$, we obtain the inequality
\[\Phi(\phi) \le \limsup_{r \to \infty} \overline{F}_r(\phi) = \overline{F}_u(\phi). \]
We complete the proof.
\end{prf}
Combining Corollary 2.1 and Theorem 3.5, we obtain immediately the following reuslt.
\begin{cor}
Let $f \in P(\mathbb{R})$ and $\phi \in L^{\infty}_R(\mathbb{R})$. Then, we have
\[\overline{F}_u(\phi) = \sup_{\Phi \in \mathcal{T}(G)} \Phi(\phi), \]
and 
\[\underline{F}_u(\phi) = \inf_{\Phi \in \mathcal{T}(G)} \Phi(\phi). \]
\end{cor}

Following result gives an analytic condition under which a given function is almost convergent.
\begin{thm}
Let $f \in P(\mathbb{R})$. For a function $\phi$ in $L^{\infty}(\mathbb{R})$, $\phi$ is almost convergent to a number $\alpha$ if and only if
\[\lim_{r \to \infty} (f_r * \phi)(x) = \alpha \]
uniformly in $x \in \mathbb{R}$. In other words, 
\[\lim_{r \to \infty} \|f_r * \phi - \alpha\|_{\infty} = 0. \]
\end{thm}

\begin{prf}
We begin with showing an elementary fact that for a real-valued essentially bounded funcition $\phi$, $\lim_{r \to \infty} (f_r*\phi)(x) = \alpha$ uniformly in $x \in \mathbb{R}$ if and only if $\overline{F}_u(\phi) = \underline{F}_u(\phi) = \alpha$, that is,
\[\limsup_{r \to \infty} \sup_{x \in \mathbb{R}} (f_r * \phi)(x) = \liminf_{r \to \infty} \inf_{x \in \mathbb{R}} (f_r * \phi)(x) = \alpha. \]
In fact, if the above equation holds, then, for any fixed $\varepsilon > 0$, we can choose a number $R$ such that 
\[\alpha - \varepsilon < \inf_{r \ge R} \inf_{x \in \mathbb{R}} (f_r * \phi)(x) \le \sup_{r \ge R} \sup_{x \in \mathbb{R}} (f_r * \phi)(x) < \alpha + \varepsilon. \]
Hence, if $r \ge R$, we have $|(f_r * \phi)(x) - \alpha| < \varepsilon$ for every $x \in \mathbb{R}$. This shows uniform convergence of $(f_r * \phi)(x)$ to $\alpha$ as $r \to \infty$.
 
On the other hand, suppose that $\lim_{r \to \infty} (f_r * \phi)(x) = \alpha$ uniformly in $x \in \mathbb{R}$. Then, for any positive number $\varepsilon > 0$, there exists a number $R > 0$ such that if $r \ge R$, 
\[|(f_r * \phi)(x) - \alpha| < \varepsilon \]
for every $x \in \mathbb{R}$. Namely, 
\[\alpha - \varepsilon < (f_r * \phi)(x) < \alpha + \varepsilon \]
for every $x \in \mathbb{R}$. This means that
\[\alpha - \varepsilon < \inf_{x \in \mathbb{R}} (f_r * \phi)(x) \le \sup_{x \in \mathbb{R}} (f_r * \phi)(x) < \alpha + \varepsilon \]
whenever $r \ge R$. Considering the limit as $r \to \infty$, we obtain the equation 
\[\liminf_{r \to \infty} \inf_{x \in \mathbb{R}} (f_r * \phi)(x) = \limsup_{r \to \infty} \sup_{x \in \mathbb{R}} (f_r * \phi)(x) = \alpha. \]
We obtain the desired result. 

Now suppose that $\phi \in L^{\infty}_R(\mathbb{R})$ is almost convergent to $\alpha$. Hence, by definition of almost convergence, we have 
\[\Phi(\phi) = \alpha \]
for every $\Phi \in \mathcal{T}$. We show that 
\[\limsup_{r \to \infty} \sup_{x \in \mathbb{R}} (f_r * \phi)(x) = \liminf_{r \to \infty} \inf_{x \in \mathbb{R}} (f_r * \phi)(x) = \alpha. \]
By Theorem 3.4, for any $\Phi \in \mathcal{T}$, it holds that 
\[\liminf_{r \to \infty} \inf_{x \in \mathbb{R}} (f_r * \phi)(x) \le \Phi(\phi) \le \limsup_{r \to \infty} \sup_{x \in \mathbb{R}} (f_r * \phi)(x). \]
Now assume that $\beta := \liminf_{r \to \infty} \inf_{x \in \mathbb{R}} (f_r * \phi)(x) < \limsup_{r \to \infty} \sup_{x \in \mathbb{R}} (f_r * \phi)(x) =: \gamma$ is true. Then, for any real number $\rho \in [\beta, \gamma]$, there exists a topologically invariant mean $\Phi$ such that $\Phi(\phi) = \rho$. In fact, define $\Phi_0$ on the subspace $\mathbb{R}\phi = \{c\phi : c \in \mathbb{R}\}$ of $L_R^{\infty}(\mathbb{R})$ by $\Phi_0(c\phi) = c\rho$. Then, it is easily confirmed that
\[\Phi_0(\psi) \le \overline{F}_u(\psi) \]
on $\mathbb{R}\phi$. By the Hahn-Banach extension theorem, $\Phi_0$ extends to a linear functional $\Phi$ on $L^{\infty}_R(\mathbb{R})$ satisfying 
\[\Phi(\psi) \le \overline{F}_u(\psi) \]
for every $\psi \in L^{\infty}_R(\mathbb{R})$. Again by Theorem 3.4, this $\Phi$ is topologically invariant. This contradicts the assumption that $\Phi(\phi) =\alpha$ for every $\Phi \in \mathcal{T}$.

Conversely, for a function $\phi \in L^{\infty}_R(\mathbb{R})$, suppose that $\lim_{r \to \infty} \|f_r*\phi-\phi\|_
{\infty} = 0$, that is, $\overline{F}_u(\phi) = \underline{F}_u(\phi) = \alpha$. Then, by Theorem 3.4, we obtain immediately $\Phi(\phi) = \alpha$ for every $\Phi$ in $\mathcal{T}$, which means that $\phi$ is almost convergent to the number $\alpha$.

For a general complex-valued function $\phi = u + iv \in L^{\infty}(\mathbb{R})$, it follows easily from the real-valued case by the fact that $\phi$ is almost covergent to the complex number $\alpha = a + bi$ if and only if $u$ and $v$ are almost convergent to the real numbers $a$ and $b$, respectively.
\end{prf}

\begin{cor}
Let $\phi \in L^{\infty}(\mathbb{R})$. Then, $\phi$ is almost convergent to a number $\alpha$ if and only if 
\[\lim_{\theta \to \infty} \frac{1}{2\theta} \int_{x-\theta}^{x+\theta} \phi(t)dt = \alpha \]
uniformly in $x \in \mathbb{R}$. Or, equivalently, 
\[\lim_{x \to \infty} \frac{1}{\pi} \int_{-\infty}^{\infty} \phi(t) \frac{x}{x^2+(y-t)^2}dt = \alpha \]
uniformly in $y \in \mathbb{R}$.
\end{cor}
We note that the first half of Corollary 3.2 was obtained by Raimi [13]. Another kind of analytic condition than Theorem 3.5 can be found in [8].

\section{An application to the bounded analytic functions on the half plane}
In this section, we consider some applications of almost convergence to the class of bounded analytic functions $H^{\infty}(\mathbb{C}^+)$ on the right half plane $\mathbb{C}^+ := \{z \in \mathbb{C} : Re (z) > 0\}$. As a reference of Hardy space, we refer the reader to [4].
Recall that each $\phi \in H^{\infty}(\mathbb{C}^+)$ has nontangential limits $\lim_{x \to 0^+} \phi(z) = \lim_{x \to 0^+} \phi(x+iy)$ almost everywhere on $\mathbb{R}$ and we can identify a function in $H^{\infty}(\mathbb{C}^+)$ with its boundary function $\phi(iy) \in L^{\infty}(\mathbb{R})$. We denote by $H^{\infty}(\mathbb{R})$ the boundary functions of $H^{\infty}(\mathbb{C}^+)$ and we have the isometry $H^{\infty}(\mathbb{C}^+) \cong H^{\infty}(\mathbb{R})$ as Banach algebras. Recall that if a function $\phi$ in $H^{\infty}(\mathbb{R})$ is given, then we can extend $\phi$ to a bounded analytic function $\hat{\phi}$ on $\mathbb{C}^+$ by the Poisson integral of $\phi$:
\[\hat{\phi}(z) = \frac{1}{\pi}\int_{-\infty}^{\infty} \phi(t) P_x(y-t)dt = \frac{1}{\pi}\int_{-\infty}^{\infty} \phi(t) \frac{x}{x^2 + (y-t)^2}dt \quad (z = x + iy, \ x > 0, \ y \in \mathbb{R}). \]

Using Corollaries $3.1$ and $3.2$, we can relate the behavior of $\phi \in H^{\infty}(\mathbb{C}^+)$ on the imaginary axis to that at infinity.
\begin{thm}
Let $\phi \in H^{\infty}(\mathbb{C}^+)$. Then, we have the following equations:
\[\limsup_{x \to \infty} \sup_{y \in \mathbb{R}} \frac{1}{2\theta} \int_{y-\theta}^{y+\theta} Re(\phi(it))dt = \limsup_{x \to \infty} \sup_{y \in \mathbb{R}} Re(\phi(x+iy)), \]
and 
\[\liminf_{x \to \infty} \inf_{y \in \mathbb{R}} \frac{1}{2\theta} \int_{y-\theta}^{y+\theta} Re(\phi(it))dt = \liminf_{x \to \infty} \inf_{y \in \mathbb{R}} Re(\phi(x+iy)), \]
We also have the same relation with respect to the imaginary part of $\phi$.
\end{thm}

\begin{thm}
Let $\phi \in H^{\infty}(\mathbb{C}^+)$. Then, the boundary function $\phi(iy)$ is almost convergent to a number $\alpha$, that is, 
\[\lim_{\theta \to \infty} \frac{1}{2\theta} \int_{y-\theta}^{y+\theta} \phi(it)dt = \alpha \]
uniformly in $y \in \mathbb{R}$ if and only if \[\lim_{x \to \infty} \phi(x+iy) = \alpha \]
uniformly in $y \in \mathbb{R}$.
\end{thm}

We now consider the interpretation of Theorem 4.1 from the Banach algebra theory. Let $\mathfrak{M}$ be the maximal ideal space of the Banach algebra $H^{\infty}(\mathbb{C}^+)$. Note that a maximal ideal $\mathfrak{m}$ of $H^{\infty}(\mathbb{C}^+)$ is equivalent to the complex homomorphism $\chi$ of $H^{\infty}(\mathbb{C}^+)$ onto $\mathbb{C}$ induced by the canonical mapping $H^{\infty}(\mathbb{C}^+) \rightarrow H^{\infty}(\mathbb{C}^+)/\mathfrak{m} \cong \mathbb{C}$. Identifying the evaluation mappings $\hat{z} : H^{\infty}(\mathbb{C}^+) \ni \phi \mapsto \phi(z)$ at the point $z \in \mathbb{C}^+$ with the point $z \in \mathbb{C}^+$ itself, we can see $\mathbb{C}^+$ as a subset of $\mathfrak{M}$. Then, it follows that $\mathbb{C}^+$ is dense in $\mathfrak{M}$ by the corona theorem. In other words, for each element $\chi \in \mathfrak{M}$, there exists a net $\{z_{\alpha}\} \in \mathbb{C}^+$ such that $\chi = \lim_{\alpha} \hat{z}_{\alpha}$ in a weak* sense. Let $\mathfrak{M}_{\infty}$ be the subset of $\mathfrak{M}$ consisting of those elements $\chi$ which are the limit of a net $\{\hat{z}_{\alpha}\}$ with $\lim_{\alpha} x_{\alpha} = \infty$, where $z_{\alpha} = x_{\alpha} + iy_{\alpha}$.

\begin{thm}
Let $\chi$ be in $\mathfrak{M}_{\infty}$. Then, there exists a topologically invariant mean $\Phi$ such that $\Phi|_{H^{\infty}(\mathbb{C}^+)} = \chi$.
\end{thm}

\begin{prf}
Let $\{z_{\alpha}\}$ be a net such that $\chi = w^*\mathchar`-\lim_{\alpha} \hat{z}_{\alpha}$ in $H^{\infty}(\mathbb{C}^+)$. Taking subnet of $\{\hat{z}_{\alpha}\}$ if necessary so that the net $\{\hat{z}_{\alpha}\}$ converges to an element $\Phi$ in $L^{\infty}(\mathbb{R})^*$ in a weak* sense. 
Note that
\[\hat{z}_{\alpha}(\phi) = \phi(z_{\alpha}) = \frac{1}{\pi}\int_{-\infty}^{\infty} \phi(t) \frac{x_{\alpha}}{x_{\alpha}^2 + (y_{\alpha}-t)^2}dt \quad (z_{\alpha} = x_{\alpha} + iy_{\alpha}), \]
where $\phi \in L^{\infty}(\mathbb{R})$. For any fixed $x > 0$ and $\phi \in L^{\infty}_R(\mathbb{R})$, we have
\[\sup_{y \in \mathbb{R}} \frac{1}{\pi} \int_{-\infty}^{\infty} \phi(t) \frac{x_{\alpha}}{x_{\alpha}^2 + (y_{\alpha}-t)^2}dt \le \sup_{y \in \mathbb{R}} \frac{1}{\pi} \int_{-\infty}^{\infty} \phi(t) \frac{x}{x^2 + (y-t)^2}dt \]
whenever $x_{\alpha} \ge x$. Hence, for each $\phi \in L^{\infty}_R(\mathbb{R})$, 
\[\Phi(\phi) = \lim_{\alpha} \frac{1}{\pi} \int_{-\infty}^{\infty} \phi(t) \frac{x_{\alpha}}{x_{\alpha}^2 + (y_{\alpha}-t)^2}dt \le \overline{P}_x(\phi) \]
is valid for any $x > 0$. Taking the limit superior as $x \to \infty$, we obtain
\[\Phi(\phi) \le \overline{P}_u(\phi), \]
which shows that $\Phi$ is in $\mathcal{T}$ by Theorem $3.4$. By the construction of $\Phi$, it is obvious that $\Phi|_{H^{\infty}(\mathbb{C}^+)} = \chi$.
\end{prf}

\begin{cor}
Almost convergence is multiplicative on $H^{\infty}(\mathbb{R})$. That is, if $\phi$ and $\psi$ in $H^{\infty}(\mathbb{R})$ almost converges to $\alpha$ and $\beta$, respectively, then $\phi\psi$ almost converges to $\alpha\beta$.
\end{cor}

In fact, in the more general situation than Theorem 4.2, the case in which the uniform limit $\lim_{x \to \infty} \phi(x+iy)$ does not exist, the relation between the values $\{\Phi(\phi)\}_{\Phi \in \mathcal{T}(\mathbb{R})}$ and the cluster points of $\{\phi(z)\}_{z \in \mathbb{C}^+}$ is valid. 
\begin{cor}
Let $\phi \in H^{\infty}(\mathbb{C}^+)$. Let $\alpha$ be a cluster point of $\phi(z)$ as $Re (z)$ tends to infinity, that is, there exists a sequence $\{z_n\}_{n=1}^{\infty} = \{x_n+iy_n\}_{n=1}^{\infty}$ of $\mathbb{C}^+$ with $\lim_{n \to \infty} x_n = \infty$ such that $\lim_{n \to \infty} \phi(x_n + iy_n) = \alpha$. Then, there exists a topologically invariant mean $\Phi \in \mathcal{T}$ with $\Phi(\phi) = \alpha$. Here, $\phi = \phi(iy)$ is the boundary function on the imaginary axis.
\end{cor}
The proof is obvious from Theorem 4.3. Theorem 4.1 and Corollary 4.2 briefly demonstrate the importance of the topologically invariant means in the study of the Hardy space $H^{\infty}(\mathbb{C}^+)$, especially, for the investigation of the behavior of bounded analytic functions at infinity. More detailed results on the maximal ideal space of the Hardy space on the half space can be found, for example, in [11], [17], [18].

Next, we discuss the relationship of the above results to Wiener's tauberian theorem. See [15] for the proof of the following theorem.
\begin{thm}[Wiener's Tauberian theorem]
Let $G$ be a locally compact abelian group and $m$ be the Haar measure of $G$. Suppose $f$ be in the group algebra $L^1(G)$ of $G$ such that its Fourier transform $\hat{f}$ does not vanish on the dual group $\hat{G}$ of $G$. Then, for any essentially bounded measurable function $\phi$ on $G$, if
\[\lim_{x \to \infty} f * \phi(x) = \int_G f(xt^{-1})\phi(t)dt = \hat{f}(0)\alpha \]
holds, then for any $g \in L^1(G)$, we have
\[\lim_{x \to \infty} g * \phi(x) = \int_G g(xt^{-1})\phi(t)dt = \hat{g}(0)\alpha. \]
\end{thm}

We now apply this theorem to the positive part $\mathbb{R}^{\times}_+ = (0, \infty)$ of the multiplicative group $\mathbb{R}^{\times}$ of $\mathbb{R}$. Note that the Haar measure of $\mathbb{R}^{\times}$ is $dt/|t|$.
Now, we determine the dual group of $\mathbb{R}^{\times}_+$, that is, the characters on $\mathbb{R}^{\times}_+$. Observe that $\mathbb{R}^{\times}_+$ is isomorphic to the additive group $\mathbb{R}$ via the isomorphism $\mathbb{R} \ni x \mapsto \log x \in \mathbb{R}^{\times}_+$. Through this isomorphism, each character $\chi$ on $\mathbb{R}$ induces a character $\chi^{\times}$ by the composition $\chi \circ \log$. 
\[
\begin{diagram}
    \node{\mathbb{R}^{\times}_+} \arrow{e,t}{\chi^{\times}} \arrow{s,l}{\log x} \node{\mathbb{C}^{\times}} \\
    \node{\mathbb{R}} \arrow{ne} 
\end{diagram}
\]
Conversely, any character on $\mathbb{R}_+^{\times}$ is obtained in this way. Hence, each character $\chi^{\times}$ on $\mathbb{R}^{\times}_+$ can be written as $e^{i\xi \log x} = x^{i\xi} \ (\xi \in \mathbb{R})$. 

\begin{cor}[Wiener's Tauberian theorem for $\mathbb{R}^{\times}_+$]
Suppose $f$ be in $L^1(\mathbb{R}^{\times}_+)$ with $\hat{f}$ does not vanish on $\hat{\mathbb{R}}^{\times}_+$. That is, 
\[\hat{f}(\xi) = \int_0^{\infty} f(t)t^{i\xi} \frac{dt}{t} \not= 0 \]
for every $\xi \in \mathbb{R}$. For $\phi \in L^{\infty}(\mathbb{R}_+^{\times})$, if 
\[\lim_{x \to \gamma} f \overset{M}{*} \phi(x) = \lim_{x \to \gamma} \int_0^{\infty} \phi(t) f\left(x/t\right)\frac{dt}{t} = \hat{f}(0)\alpha  \]
holds, then for every $g \in L^1(\mathbb{R}^{\times}_+)$, we have
\[\lim_{x \to \gamma} g \overset{M}{*} \phi(x) = \lim_{x \to \gamma} \int_0^{\infty} \phi(t) g\left(x/t\right)\frac{dt}{t} = \hat{g}(0)\alpha.  \]
Here, $\gamma = 0^+$ or $\infty$ and the symbol $\overset{M}{*}$ stands for convolution in $L^1(\mathbb{R}^{\times})$ sometimes called Mellin convolution.
\end{cor}

\bigskip
Observe that, for any $r > 0$, $f \in P(\mathbb{R})$ and $\phi \in L^{\infty}(\mathbb{R})$, the convolution $f_r * \phi(x)$ can be transformed as follows:
\begin{align}
f_r * \phi(x) &= \int_{-\infty}^{\infty} \phi(x-t)f_r(t)dt \notag \\
&= \int_{-\infty}^{\infty} \phi(x-t)\frac{1}{r} f\left(\frac{t}{r}\right)dt \notag \\
&= \int_{-\infty}^{\infty} \phi(x-t)\frac{|t|}{r} f\left(\frac{t}{r}\right)\frac{dt}{|t|} \notag \\
&= \int_{-\infty}^{\infty} \phi(x-t) \tilde{f}\left(\frac{r}{t}\right)\frac{dt}{|t|} \quad (\tilde{f}(t) := |t|^{-1}f(t^{-1})) \notag \\
&= \int_{-\infty}^{\infty} \phi_x(-t) \tilde{f}\left(\frac{r}{t}\right)\frac{dt}{|t|} \notag \\
&= (\tilde{f} \overset{M}{*} {\phi_x}^*)(r) \quad (\phi^*(t) := \phi(-t)). \notag
\end{align}
Here note that $\tilde{f}$ is in $L^1(\mathbb{R}^{\times})$. Additionally, if we assume $f \in P(\mathbb{R})$ is even, that is, $f(-x) = f(x)$, the above formula can be written as
\begin{align}
f_r * \phi(x) &= \int_0^{\infty} \{{\phi_x}^*(t) + {\phi_x}^*(-t)\}\tilde{f}\left(\frac{r}{t}\right) \frac{dt}{t} \notag \\ 
&= \int_{\mathbb{R}^{\times}_+} {\phi_x}^{\sharp}(t)\tilde{f}\left(\frac{r}{t}\right) \frac{dt}{t} \quad (\phi^{\sharp}(t) := \phi^*(t) + \phi^*(-t)) \notag \\
&= (\tilde{f} \overset{M}{*} {\phi_x}^{\sharp})(r). \notag
\end{align}

Combining  Corollary 4.3 with the above observation, we obtain the following result.
\begin{thm}
Let $\phi \in L^{\infty}(\mathbb{R})$ and $f \in P(\mathbb{R})$ be an even function such that 
\[\int_0^{\infty} f(t)t^{i\xi}dt \]
does not vanish for every $\xi \in \mathbb{R}$. Then, if for any fixed $x \in \mathbb{R}$, 
\[\lim_{r \to \gamma} f_r * \phi(x) = \alpha \]
holds, then for every even function $g \in P(\mathbb{R})$, we have
\[\lim_{r \to \gamma} g_r * \phi(x) = \alpha, \]
where $\gamma = 0^+$ or $\infty$.
\end{thm}

\begin{prf}
Suppose that $f$ be in $P(\mathbb{R})$ satisfying the condition of the theorem. Observe that 
\begin{align}
\int_0^{\infty} \tilde{f}(t)t^{i\xi}\frac{dt}{t} &= \int_0^{\infty} t^{-1}f(t^{-1})t^{i\xi}\frac{dt}{t} \notag \\
&= \int_0^{\infty} tf(t)t^{-i\xi}\frac{dt}{t}\notag \\
&= \int_0^{\infty} f(t)t^{-i\xi}dt \notag
\end{align}
and we have $\hat{\tilde{f}}(\xi) \not= 0$ for every $\xi \in \mathbb{R}$ by the assumption on $f$. Hence, applying Corollary $4.3$,  for any $\phi \in L^{\infty}(\mathbb{R})$, we conclude that if 
\[\lim_{r \to \infty} f_r * \phi(x) = \lim_{r \to \infty} (\tilde{f} \overset{M}{*} \phi_x^{\sharp})(r) = \alpha \]
holds, then we have
\[\lim_{r \to \infty} g_r * \phi(x) = \lim_{r \to \infty} (\tilde{g} \overset{M}{*} \phi_x^{\sharp})(r) = \alpha. \]
This completes the proof.
\end{prf}

Now consider the cases $f(x) = D(x) = \frac{1}{2}I_{[-1,1]}(x)$ and $f(x) = P(x) = \frac{1}{\pi} \frac{1}{1+x^2}$. First, we  calculate the following integral:
\[\int_0^{\infty} f(t)t^{i\xi}dt. \]
For the function $D(x)$, we have
\begin{align}
\frac{1}{2} \int_0^{\infty} I_{[-1,1]}(t)t^{i\xi}dt &= \frac{1}{2} \int_0^1 t^{i\xi}dt \notag \\
&= \frac{1}{2} \left[\frac{t^{i\xi+1}}{i\xi+1}\right]_0^1 \notag \\
&= \frac{1}{2} \frac{1}{i\xi+1}, \notag
\end{align}
which does not vanish for every $\xi \in \mathbb{R}$. For the function $P(x)$, by integration by substitution, we obtain
\[\int_0^{\infty} \frac{1}{\pi}\frac{1}{1+t^2} t^{i\xi}dt = \frac{1}{2\pi} \Gamma\left(\frac{1-i(\xi-1)}{2}\right)\Gamma\left(\frac{1+i(\xi+1)}{2}\right), \]
where $\Gamma$ is as usual the gamma function. Thus, we see that the integral does not vanish for every $\xi \in \mathbb{R}$.
As a result, the following is obtained from Theorem 4.4.
\begin{thm}
Let $\phi \in H^{\infty}(\mathbb{C}^+)$. Then, the following three conditions are equivalent. \\
$(1)$ $\lim_{\theta \to \infty} \frac{1}{2\theta} \int_{-\theta}^{\theta} \phi(it)dt = \alpha$. \\
$(2)$ $\lim_{x \to \infty} \phi(x+iy) = \alpha$ for some $y \in \mathbb{R}$. \\
$(3)$ $\lim_{x \to \infty} \phi(x+iy) = \alpha$ for every $y \in \mathbb{R}$.
\end{thm}

\begin{prf}
$(1) \Leftrightarrow (3)$. For each $y \in \mathbb{R}$, observe that 
\[D_{\theta} * \phi(y) = \frac{1}{2\theta} \int_{y-\theta}^{y+\theta} \phi(t)dt \ (\theta > 0), \quad P_x * \phi(y) = \frac{1}{\pi} \int_{-\infty}^{\infty} \phi(t)\frac{x}{x^2+ (y-t)^2}dt \ (x > 0)  \]
and by Theorem $4.5$, it follows that
\[\lim_{\theta \to \infty} \frac{1}{2\theta} \int_{y-\theta}^{y+\theta} \phi(t)dt = \alpha \]
if and only if 
\[\lim_{x \to \infty} \phi(x+iy) = \alpha \]
Furthermore, it is obvious that
\[\lim_{\theta \to \infty} D_{\theta} * \phi(y) = \alpha \Leftrightarrow \lim_{\theta \to \infty} D_{\theta} * \phi(0) = \alpha. \]
We thus obtain $(1) \Leftrightarrow (3)$ and $(2) \Rightarrow (1)$. The implication $(3) \Rightarrow (2)$ is obvious.  We complete the proof.
\end{prf}

\begin{cor}
The continuous version of Ces\`{a}ro mean is multiplicative on $H^{\infty}(\mathbb{R})$. That is, if functions $\phi, \psi \in H^{\infty}(\mathbb{R})$ has the following limits
\[\lim_{\theta \to \infty} \frac{1}{2\theta} \int_{-\theta}^{\theta} \phi(t)dt = \alpha, \quad \lim_{\theta \to \infty} \frac{1}{2\theta} \int_{-\theta}^{\theta} \psi(t)dt = \beta, \]
then the product $\phi\psi$ satisfies the following equation
\[\lim_{\theta \to \infty} \frac{1}{2\theta} \int_{-\theta}^{\theta} \phi(t)\psi(t)dt = \alpha\beta. \]
\end{cor}

Finally, we close with the remark on the relation to Fatou's theorem. Note that, as was stated in [6], for each element $f \in P(\mathbb{R})$, the functions $\{f_r\}_{r > 0}$ form an approximate identity of $L^1(\mathbb{R})$. By the Lebesgue differentiation theorem,  for any $\phi \in L^{\infty}(\mathbb{R})$, we obtain 
\[\lim_{\theta \to 0} (D_{\theta} * \phi)(x) = \lim_{\theta \to 0} \frac{1}{2\theta} \int_{x-\theta}^{x+\theta} \phi(t)dt = \phi(x) \]
exists a.e. on $\mathbb{R}$. By theorem 4.5, we immeidiately obtain the following reuslt:
\begin{cor}
Let $f \in P(\mathbb{R})$ be an even function and consider the approximate identity $\{f_r\}_{r > 0}$. Then, for any $\phi \in L^{\infty}(\mathbb{R})$, we have
\[\lim_{r \to 0^+} (f_r * \phi)(x) = \phi(x) \]
a.e. on $\mathbb{R}$. 
\end{cor}

In particular, we obtain the special case of Fatou's theorem, in which an integrand is a function in $L^{\infty}(\mathbb{R})$. See [7], [10], [15] for more comprehensive arguments about the inverse of Fatou's theorem.
\begin{thm}
Let $\phi \in L^{\infty}(\mathbb{R})$. We define the harmonic function $\hat{\phi}(z)$ on $\mathbb{C}^+$ by the Poisson integral.
Then, for any $y \in \mathbb{R}$, the radial limit
\[\lim_{x \to 0+} \hat{\phi}(x+iy) \]
exists if and only if the symmetric difference 
\[\lim_{\theta \to 0} \frac{1}{2\theta} \int_{x-\theta}^{x+\theta} \phi(t)dt \]
exists.
\end{thm}

\end{document}